\documentclass[11pt]{amsart}
\usepackage{mesfonts}

\def\crit{{\mathrm{Crit}}}

\def\argmin{\mathop{\mathrm{argmin}}}

\RequirePackage{amssymb}

\newtheorem{theorem}{Theorem}

\newtheorem{corollary}{Corollary}

\newtheorem{assumption}{Assumption}

\usepackage{graphicx}
\RequirePackage{amssymb}

\begin{document}

\title{A pseudo-RIP for multivariate regression}
\author{Christophe Giraud}

\date{June 2011}
\begin{abstract}
We give a suitable RI-Property under which recent results for trace regression translate into strong risk bounds for multivariate regression. This pseudo-RIP is compatible with the setting $n<p$.
\end{abstract}
\keywords{Multivariate regression, Restricted Isometry Property\\ \\ {\small \emph{Institute: }Ecole Polytechnique, CMAP, UMR CNRS 7641 (FRANCE)} \vspace{0.1cm}
\\ {\small \emph{E-mail address: }christophe.giraud@polytechnique.edu}}

\maketitle

\section{Introduction}
\subsection{Statistical framework}
Multivariate regression deals with $n$ observations of a $T$-dimensional vector
$$y_{i}=A_{0}^Tx_{i}+\varepsilon_{i},\quad i=1,\ldots,n$$
where $A_{0}^T$ is the transpose of a $p\times T$ matrix $A_{0}$. We have in mind that $A_{0}$ has a small (unknown) rank and the design $x_{i}$ is non-random.
Writing $Y$, $X$ and $E$ for the matrices with respective rows $y_{i}^T$, $x_{i}^T$ and $\varepsilon_{i}^T$, the above equation translate into
$$Y=XA_{0}+E.$$
Anderson~\cite{Anderson51} and Izenman~\cite{Izenman75} have introduced reduced-rank estimators
\begin{equation*}\label{RR}
\widehat A_{r}\in\argmin_{A\ :\ \mathrm{rank}(A)\leq r} \|Y-XA\|^2,\quad r=0,\ldots,\min(p,T),
\end{equation*}
where $\|.\|$ is the Hilbert-Schmidt norm associated to the scalar product $\langle .,.\rangle$. The problem of selecting among the family of estimators $\ac{\widehat A_{r},\ r=0,\ldots,\min(p,T)}$ by minimizing the criterion
\begin{equation*}\label{critAr}
\crit(r)=\|Y-X\widehat A_{r}\|^2+\pen(r)\sigma^2\quad\textrm{and}\quad \crit'(r)=\log\pa{\|Y-X\widehat A_{r}\|^2}+\pen'(r)
\end{equation*}
has been investigated recently from a non-asymptotic point of view by Bunea {\it et al.}~\cite{BSW11} and Giraud~\cite{Giraud10}. Both papers provide oracle bounds for the predictive risk
$$R(\widehat A)=\E\cro{\|X\widehat A-XA\|^2}$$
with \emph{no assumption} on the design $X$.
\bigskip

Multivariate regression corresponds to a special case of the trace regression model
$$\mathcal{Y}_{j}=\langle Z_{j},A_{0}\rangle+\xi_{j},\quad j=1,\ldots,N,$$
where $\langle Z_{j},A_{0}\rangle=\mathrm{tr}(Z_{j}^TA_{0})$. Indeed, we have for all $i\in\ac{1,\ldots,n}$ and $t\in \ac{1,\ldots,T}$
\begin{eqnarray*}
Y_{it}&=&\langle A_{0}^Tx_{i},e_{t} \rangle + E_{it}\\
&=& \langle\,\underbrace{x_{i}e_{t}^T}_{=:Z_{it}}\,,A_{0} \rangle + E_{it},
\end{eqnarray*}
where $\ac{e_{1},\ldots,e_{T}}$ is the canonical basis of $\R^T$. Many recent works~\cite{Bach08,NeghabanWainwright,
 RohdeTsybakov, KLT11} have investigated trace regression with nuclear norm penalization. Translated in terms of multivariate regression, Nuclear-Norm-Penalized regression estimators are defined by
\begin{equation}\label{NNP}
\widehat A_{\lambda}\in\argmin_{A\in\R^{p\times T}} \ac{\|Y-XA\|^2+\lambda \sum_{k}\sigma_{k}(A)},
\end{equation}
where $\sigma_{1}(A)\geq\sigma_{2}(A)\geq \ldots$ are the singular values of $A$. Several risk bounds have been obtained for the predictive risk of $\widehat A_{\lambda}$ and they all require the assumption (semi-RI Property)
\begin{equation}\label{semiRIP}
\|A\|\leq \mu\, \|XA\|,\quad \textrm{for all}\ A\in\R^{p\times T} 
\end{equation}
for some positive $\mu$. In other words the smallest eigenvalue of $X^TX$ must be larger than $1/ \mu^2>0$. This enforces the sample size $n$ to be larger than the number $p$ of parameters. This assumption on the design needed for $\widehat A_{\lambda}$ is thus very strong, in contrast with the reduced-rank estimator $\widehat A_{\hat r}$ which requires no assumption on the design.

\subsection{Object of this note}
In this note, we emphasize that the Assumption~\eref{semiRIP} coming from the general trace regression framework can be weaken for the multivariate regression framework. Under this (much) weaker assumption, we show that the analysis  of Kolchinskii {\it et al.}~\cite{KLT11} gives an  oracle bound with leading constant 1 for the estimators $\widehat A_{\lambda}$.

\section{Semi-RIP for multivariate regression}
The Condition~\eref{semiRIP} requires the sample size $n$ to be larger than the number $p$ of covariates. 
Is-it still possible to get an oracle bound on $\|X\widehat A_{\lambda}-XA\|^2$ when $n$ is smaller than $p$ ?
\medskip

The analysis of Theorem 12 in Bunea {\it et al.}~\cite{BSW11} suggests that the Condition~\eref{semiRIP} only need to hold true for matrices $A$ of rank at most twice the rank of $A_{0}$. Unfortunately, when the rank of $A_{0}$ is positive this condition is still equivalent to require that  the smallest eigenvalue of $X^TX$ is larger than $1/ \mu^2>0$.
\medskip

In the analysis of in  Kolchinskii {\it et al.}~\cite{KLT11},
the Condition~\eref{semiRIP} is needed for comparing $\|\widehat A_{\lambda}-A\|$  to $\|X\widehat A_{\lambda}-XA\|$, see for example the Display~(2.17) of~\cite{KLT11}. We point out below, that this inequality needs not to hold for all matrices $\widehat A_{\lambda}$ and $A$, so that  Condition~\eref{semiRIP} can be relaxed to handle cases where $p>n$. 
\begin{assumption}
$$\sigma_{\mathrm{rank}(X)}(X)\geq {1\over \mu}>0$$
where $\sigma_{1}(X)\geq \sigma_{2}(X)\geq \ldots$ are the singular values of $X$. 
\end{assumption}
The singular value $\sigma_{\mathrm{rank}(X)}(X)$ is always positive but can be arbitrary small. Assumption~1 requires a positive lower bound on this singular value.

\subsection{Risk bound under Assumption 1} 
Write $\mathrm{rg}(X^T)$ for the range of the linear operator $X^T$ and $\Pi_{\mathrm{rg}(X^T)}$ for the orthogonal projection onto the range of $X^T$ in $\R^p$.
Since we have the orthogonal decomposition $\R^p=\mathrm{ker}(X)+\mathrm{rg}(X^T)$, we have
$X\Pi_{\mathrm{rg}(X^T)}A=XA$ for any matrix $A$. In addition, $\sigma_{k}(\Pi_{\mathrm{rg}(X^T)}A)\leq \sigma_{k}(A)$ for any $k$ and matrix $A$, so 
$$\sum_{k} \sigma_{k}(\Pi_{\mathrm{rg}(X^T)}A) \leq  \sum_{k} \sigma_{k}(A),$$
with strict inequality if $\Pi_{\mathrm{rg}(X^T)}A\neq A$. As a consequence, we have $\Pi_{\mathrm{rg}(X^T)}\widehat A_{\lambda}=\widehat A_{\lambda}$, so
 $\widehat A_{\lambda}$ is also a minimizer of
\begin{equation}\label{NNP2}
\min_{A\in\mathbb{A}} \ac{\|Y-XA\|^2+\lambda \sum_{k}\sigma_{k}(A)}.
\end{equation}
where
$\mathbb{A}:=\ac{A\in\R^{p\times T}\ :\ \mathrm{rg}(A)\subset \mathrm{rg}(X^T)}.$
\bigskip

Under Assumption 1, we have 
$$\|A\|\leq \mu\, \|XA\|,\quad \textrm{for all}\ A\in\mathbb{A}.$$
Theorem 1 of~Kolchinskii {\it et al.}~\cite{KLT11} then gives the upper bound
$$\|X\widehat A_{\lambda}-XA_{0}\|^2\leq \inf_{A\in\mathbb{A}}\ac{\|XA-XA_{0}\|^2+\pa{1+\sqrt{2}\over 2}^2\mu^2\lambda^2\mathrm{rank}(A)}$$
for $\lambda\geq 2 \sigma_{1}(X^TE)$.
Again, since $X\Pi_{\mathrm{rg}(X^T)}A=XA$ and $\mathrm{rank}(\Pi_{\mathrm{rg}(X^T)}A)\leq \mathrm{rank}(A)$, the infimum on the right hand side coincides with the infimum on the whole space $\R^{p\times T}$. We then have the following result.

\begin{theorem}
Let $\widehat A_{\lambda}$ be defined by~\eref{NNP}. Then,  
 under Assumption 1, for $\lambda\geq 2 \sigma_{1}(X^TE)$ we have
\begin{eqnarray*}
\|X\widehat A_{\lambda}-XA_{0}\|^2&\leq& \inf_{A\in\R^{p\times T}}\ac{\|XA-XA_{0}\|^2+{3\over 2}\,\mu^2\lambda^2\mathrm{rank}(A)}\\
&=& \inf_{r}\bigg\{\sum_{k\geq r+1}\sigma_{k}(XA_{0})^2+{3\over 2}\,\mu^2\lambda^2r\bigg\}.
\end{eqnarray*}
\end{theorem}

\subsection{Case of Gaussian errors}
The above statement is purely deterministic.
In the case of Gaussian errors we have the following corollary.
\begin{corollary}
Assume that the entries of $E$ are i.i.d.\ with Gaussian $\mathcal{N}(0,\sigma^2)$ distribution. Let $K>1$ and set 
$$\lambda=2K\sigma_{1}(X)\pa{\sqrt{T}+\sqrt{q}}\sigma,\quad \textrm{with}\ q=\mathrm{rank}(X).$$
Then, with probability larger than $1-e^{-(K-1)^2(T+q)/2}$ we have
\begin{eqnarray}
\|X\widehat A_{\lambda}-XA_{0}\|^2&\leq& \inf_{A\in\R^{p\times T}}\ac{\|XA-XA_{0}\|^2+6K^2\,{\sigma_{1}(X)^2 \over \sigma_{q}(X)^2}\pa{\sqrt{T}+\sqrt{q}}^2\sigma^2\mathrm{rank}(A)}\nonumber\\
&=&  \inf_{r}\bigg\{\sum_{k\geq r+1}\sigma_{k}(XA_{0})^2+6K^2\,{\sigma_{1}(X)^2 \over \sigma_{q}(X)^2}\pa{\sqrt{T}+\sqrt{q}}^2\sigma^2\,r\bigg\}\label{oracle}
\end{eqnarray}
\end{corollary}

\section{Discussion}
The Assumption 1, which requires that the smallest \emph{positive} singular value of $X$ is lower bounded, 
is much weaker than the Assumption~\eref{semiRIP}. In particular, this condition is fully compatible with the setting where the sample size $n$ is smaller than the number $p$ of covariables.

 The inequality~\eref{oracle} for the Nuclear Norm Penalized estimator suggests that the suitable "RI-Property"  for prediction in multivariate regression is\medskip

{\bf RI-Property :} \ \emph{There exists $\eta\in [1,+\infty[$ such that}
\begin{equation*}
1\leq {\sigma_{1}(X)\over \sigma_{q}(X)}\leq \eta,\quad \textrm{with }\ q=\mathrm{rank}(X).
\end{equation*}

\

When this condition is met with $\eta$ of reasonable size,
 the NNP-estimator achieves under the assumptions  of Corollary~1, the oracle inequality
$$\|X\widehat A_{\lambda}-XA_{0}\|^2\leq \inf_{A\in\R^{p\times T}}\ac{\|XA-XA_{0}\|^2+6K^2\,\eta^2\pa{\sqrt{T}+\sqrt{q}}^2\sigma^2\,\mathrm{rank}(A)}$$
with probability larger than $1-e^{-(K-1)^2(T+q)/2}$.
This inequality ensures that  the NNP-estimator is adaptive rate-minimax.

\end{document}